\documentclass[12pt,a4paper]{article}
\usepackage[cp1251]{inputenc} 
\usepackage{amsfonts}
\usepackage{amsmath}

\newcommand{\di}{\displaystyle}
\newcommand{\de}{\delta}

\newcommand{\la}{\lambda}

\newcommand{\iy}{\infty}

\begin{document}

\thispagestyle{empty}
\begin{center}
{\large\bf On recovering Sturm-Liouville differential operators with
deviating argument}
\end{center}

\begin{center}
{\bf N. Bondarenko and V. Yurko}
\end{center}

\noindent {\bf Abstract. } We consider second-order functional differential
operators with a constant delay. Properties of their spectral characteristics
are obtained and a nonlinear inverse problem is studied, which consists in
recovering the operators from their spectra. We establish the uniqueness and
develop a constructive algorithm for solution of the inverse problem.

\medskip
\noindent {\bf Keywords:} differential operators, deviating argument,
nonlinear inverse spectral problem.

\medskip
\noindent {\bf AMS Mathematics Subject Classification (2010):}
34A55 34K10 34K29 47E05 34B10 34L40\\

{\bf 1. Introduction} \\

In various real-world processes, the future behavior of the system depends
not only on its present state and rate of change of the state (corresponding
to the values of the function and its derivatives at the current point),
but also on the states in the past. Such processes are described by
functional differential equations with delay, which arise in physics,
biology and especially in engineering and control theory
(see the monographs [1-2] and the references therein).

This paper concerns an inverse spectral problem for the Sturm-Liouville
equation with a constant delay. Inverse problems of spectral analysis
consist in recovering operators from their spectral characteristics.
Basic results of inverse problem theory for the classical differential
Sturm-Liouville equation are contained in the monographs [3-6]. However,
operators with a constant delay appear to be more difficult for
investigation, and inverse problems for them are studied only in
a few special cases (see [7-11]). The standard methods of inverse
problem theory, applicable to differential operators (transformation
operator method, method of spectral mappings, etc.), do not work for
operators with delay, so one needs new approaches to construct a general
spectral theory for the latter ones.

In this paper, we consider the boundary value problems $L_j(q),$
$j=0,1,$ of the form
$$
-y''(x)+q(x)y(x-a)=\la y(x),\quad x\in(0,\pi),                      \eqno(1)
$$
$$
y(0)=y^{(j)}(\pi)=0,                                                \eqno(2)
$$
where $\la$ is the spectral parameter, $a\in (0,\pi),$ $q(x)$ is
a complex-valued function such that $q(x)\in L(a,\pi)$, and
$q(x)\equiv 0$ for $x\in [0,a].$ We study the inverse spectral problem
of recovering the potential $q(x)$ from the spectra of $L_j(q).$ More
precisely, let $\{\la_{nj}\}_{n\ge 1}$, $j=0,1,$ be the eigenvalues
of $L_j(q).$ The inverse problem is formulated as follows.

\smallskip
{\bf Inverse problem 1. } Given $\{\la_{nj}\}_{n\ge 1}$, $j=0,1,$
construct $q(x).$

\smallskip
Note that in the case of ``large delay'' when $a\ge\pi/2,$ the
characteristic functions of $L_j(q)$ depend linearly on the potential
$q(x),$ i.e. the inverse problem becomes linear. This linear case
was studied in [9] and [11]. For $a<\pi/2$ the characteristic functions
depend nonlinearly on the potential, i.e. the inverse problem becomes
nonlinear. This nonlinear case is essentially more difficult for
investigating and constructing the global solution of the inverse
problem. In this paper we consider the nonlinear case.
For definiteness, let $a\in[2\pi/5,\pi/2).$ The case $a<2\pi/5$
requires a separate investigation.
In this paper we obtain a global constructive procedure for the solution
of the inverse problem and establish its uniqueness. The main results
of the paper are Theorem 1 and Algorithm 1 (see Section 3 below).\\

{\bf 2. Preliminaries} \\

In this section, we study spectral properties of the boundary value
problems (1)-(2).

Let $Y(x,\la)$ be the solution of Eq. (1) satisfying the initial conditions
$$
Y(0,\la)=0,\; Y'(0,\la)=1.
$$
The eigenvalues $\{\la_{nj}\}_{n\ge 1}$ of the boundary value problem
$L_j$ coincide with the zeros of its characteristic function
$$
\Delta_j(\la):=Y^{(j)}(\pi,\la),\quad j=0,1.
$$
The following two propositions were proved in [7].

\smallskip
{\bf Lemma 1. }{\it Fix $j=0,1.$ The boundary value problem $L_j$
has a countable set of eigenvalues $\{\la_{nj}\}_{n\ge 1}$
(counting with multiplicities), and for $n\to\iy:$}
$$
\sqrt{\la_{n0}}=n+\frac{A_0\cos na}{2\pi n}
+o\left(\frac{1}{n}\right),                                        \eqno(3)
$$
$$
\sqrt{\la_{n1}}= \left(n-\frac{1}{2}\right)+
\frac{A_0\cos(n-1/2)a}{2\pi n}+o\left(\frac{1}{n}\right),          \eqno(4)
$$
where $A_0=\di\int_a^\pi q(t)\,dt.$

{\bf Lemma 2. }{\it The specification of the spectrum
$\{\la_{nj}\}_{n\ge 0}$ uniquely determines the
characteristic function $\Delta_j(\la)$ via}
$$
\Delta_0(\la)=\pi\prod_{n=1}^\iy\frac{\la_{n0}-\la}{n^2},\quad
\Delta_1(\la)=\prod_{n=1}^\iy \frac{\la_{n1}-\la}{(n-1/2)^2}.      \eqno(5)
$$

Let us investigate the connection between the characteristic functions
$\Delta_j(\la)$ and the potential $q(x)$. Let $\la=\rho^2.$ The function
$Y(x,\la)$ satisfies the integral equation
$$
Y(x,\la)=\frac{\sin\rho x}{\rho}+\int_{a}^x g(x,t,\la)Y(t-a,\la)\,dt, \eqno(6)
$$
where $g(x,t,\la)=\di\frac{\sin\rho(x-t)}{\rho}\,q(t).$
Solving Eq.(6) we get for $x\ge 2a$:
$$
Y(x,\la)=Y_0(x,\la)+Y_1(x,\la)+Y_2(x,\la),
$$
where
$$
Y_0(x,\la)=\frac{\sin\rho x}{\rho},
$$
$$
Y_1(x,\la)=-\frac{\cos\rho(x-a)}{2\rho^2}\int_a^x q(t)\,dt
+\frac{1}{2\rho^2}\int_a^x q(t)\cos\rho(x-2t+a)\,dt,                \eqno(7)
$$
$$
Y_2(x,\la)=
\int_{2a}^x \frac{\sin\rho(x-t)}{\rho}q(t)Y_1(t-a,\la)\,dt,         \eqno(8)
$$
and consequently,
$$
\Delta_0(\la)=\frac{\sin\rho\pi}{\rho}-A_0\frac{\cos\rho(\pi-a)}{2\rho^2}
+\frac{1}{2\rho^2}\int_a^\pi q(t)\cos\rho(2t-\pi-a)\,dt+Y_2(\pi,\la),
$$
$$
\Delta_1(\la)=\cos\rho\pi+A_0\frac{\sin\rho(\pi-a)}{2\rho}
+\frac{1}{2\rho}\int_a^\pi q(t)\sin\rho(2t-\pi-a)\,dt+Y'_2(\pi,\la).
$$
Denote
$$
\Delta^{*}_0(\rho):=2\rho^2\left(\Delta_0(\la)-
\frac{\sin\rho\pi}{\rho}+A_0\frac{\cos\rho(\pi-a)}{2\rho^2}\right), \eqno(9)
$$
$$
\Delta^{*}_1(\rho):=2\rho\left(\Delta_1(\la)-
\cos\rho\pi-A_0\frac{\sin\rho(\pi-a)}{2\rho}\right),                \eqno(10)
$$
Then
$$
\Delta^{*}_0(\rho)=
\int_a^\pi q(t)\cos\rho(2t-\pi-a)\,dt+\de_0(\rho),                  \eqno(11)
$$
$$
\Delta^{*}_1(\rho)=
\int_a^\pi q(t)\sin\rho(2t-\pi-a)\,dt+\de_1(\rho)                   \eqno(12)
$$
where $\de_0(\rho)=2\rho^2 Y_2(\pi,\la),$
$\de_1(\rho)=2\rho Y'_2(\pi,\la).$ Using (7) and (8) we calculate
$$
2\rho\de_0(\rho)=-A\sin\rho(\pi-2a)
+\frac{1}{2}\int_{-(\pi-2a)}^{(\pi-2a)} Q(\xi)\sin\rho\xi\,d\xi,   \eqno(13)
$$
$$
2\rho\de_1(\rho)=-A\cos\rho(\pi-2a)
-\frac{1}{2}\int_{-(\pi-2a)}^{(\pi-2a)} Q(\xi)\cos\rho\xi\,d\xi,   \eqno(14)
$$
where
$$
A=\int_{2a}^\pi q(t)dt\int_a^{t-a} q(s)ds,\;
Q(\xi)=Q_1(\xi/2+\pi/2+a)-Q_2(\xi/2+\pi/2)-Q_3(\xi/2+\pi/2),
$$
$$
Q_1(x)=q(x)\int_a^{x-a} q(s)ds,\;
Q_2(x)=q(x)\int_{x+a}^\pi q(s)ds,\;
Q_3(x)=\int_{x+a}^\pi q(s)q(s-x)ds.
$$
Note that the characteristic functions $\Delta_j(\la),$ $j = 0,1,$
depend nonlinearly on the potential $q(x).$

In order to simplify calculations, we assume that $q(x)\in AC[a,\pi].$
The general case requires small technical modifications. Denote
$q_1(x):=q'(x).$ Taking (11)-(14) into account and applying
integration by parts, we infer
$$
4\rho\Delta^{*}_0(\rho)=B_1\sin\rho(\pi-a)-2A\sin\rho(\pi-2a)
-\int_{-(\pi-a)}^{(\pi-a)} q_0(\xi)\sin\rho\xi\,d\xi
+\int_{-(\pi-2a)}^{(\pi-2a)} Q(\xi)\sin\rho\xi\,d\xi,             \eqno(15)
$$
$$
4\rho\Delta^{*}_1(\rho)=B_2\cos\rho(\pi-a)-2A\cos\rho(\pi-2a)
+\int_{-(\pi-a)}^{(\pi-a)} q_0(\xi)\cos\rho\xi\,d\xi
-\int_{-(\pi-2a)}^{(\pi-2a)} Q(\xi)\cos\rho\xi\,d\xi,             \eqno(16)
$$
where $B_1=2(q(a)+q(\pi)),$ $B_2=2(q(a)-q(\pi)),$
$q_0(\xi)=q_1(\xi/2+\pi/2+a/2).$ Denote
$$
d_0(\rho)=4\rho\Delta^{*}_0(\rho)
-B_1\sin\rho(\pi-a)+2A\sin\rho(\pi-2a),                            \eqno(17)
$$
$$
d_1(\rho)=4\rho\Delta^{*}_1(\rho)
-B_2\cos\rho(\pi-a)+2A\cos\rho(\pi-2a).                            \eqno(18)
$$
It follows from (15)-(16) and (17)-(18) that
$$
d_0(\rho)=-\int_{-(\pi-a)}^{(\pi-a)} R(\xi)\sin\rho\xi\,d\xi,
\quad d_1(\rho)=\int_{-(\pi-a)}^{(\pi-a)} R(\xi)\cos\rho\xi\,d\xi, \eqno(19)
$$
where $R(\xi)=q_0(\xi)-Q(\xi),$ and $Q(\xi)\equiv 0$ outside the
interval $(-(\pi-2a), (\pi-2a)).$ In particular, this yields
$$
q_1(x)=R(2x-\pi-a)+Q_1(x+a/2)-Q_2(x-a/2)-Q_3(x-a/2),
\; x\in(3a/2,\pi-a/2).                                            \eqno(20)
$$

\bigskip
{\bf 3. Inverse problem} \\

In this section, we provide our main results: a constructive
procedure for solving Inverse problem~1 and the corresponding
uniqueness theorem. The solution of Inverse problem 1 can be
constructed by the following algorithm.

\smallskip
{\bf Algorithm 1.} Let the spectra $\{\la_{nj}\}_{n\ge 1},$ $j=0,1,$
be given. We then

1) Construct the characteristic functions $\Delta_j(\la),$ $j=0,1,$ by (5).

2) Find the constant $A_0,$ using (3) or (4).

3) Calculate $\Delta^{*}_j(\la),$ $j=0,1,$ according to (9) and (10).

4) Find the constants $A, \, B_1$ and $B_2,$ using (15) and (16),
and calculate $q(a) = (B_1 + B_2)/4,$ $q(\pi) = (B_1 - B_2)/4$.

5) Construct $d_0(\rho)$ and $d_1(\rho)$ via (17) and (18).

6) Calculate $R(\xi),$ using (19).

7) Find the function $q_0(\xi)$ for
$\xi\in(-(\pi-a),-(\pi-2a))\cup(\pi-2a,\pi-a),$ by $q_0(\xi)=R(\xi).$

8) Calculate $q_1(x)=q_0(2x-\pi-a)$ for $x\in(a,3a/2)\cup(\pi-a/2,\pi).$

9) Calculate
$$
   q(x) = q(a) + \int_a^x q_1(t) \, dt, \: x \in (a, 3a/2), \quad
   q(x) = q(\pi) - \int_x^{\pi} q_1(t) \, dt, \: x \in (\pi-a/2, \pi).
$$

10) Using (20) and knowledge of $q(x)$ for $x\in(a,3a/2)\cup(\pi-a/2,\pi),$
construct the function $q_1(x)$ for $x\in(3a/2,\pi-a/2)$ via
$$
q_1(x)=R(2x-\pi-a)+q(x+a/2)\int_a^{x-a/2} q(s)\,ds
$$
$$
-q(x-a/2)\int_{x+a/2}^{\pi} q(s)\,ds-\int_{x+a/2}^{\pi} q(s)q(s-x+a/2)\,ds.
$$

11) Construct $q(x)$ for $x\in(3a/2,\pi-a/2).$

\medskip
Thus, we have proved the following theorem.

\medskip
{\bf Theorem 1. }{\it The specification of two spectra
$\{\la_{nj}\}_{n\ge 1},$ $j=0,1,$ uniquely determines the potential $q(x).$
The solution of Inverse problem 1 can be found by Algorithm 1.}

\bigskip
{\bf Acknowledgments.} This research was supported in part by Russian
Foundation for Basic Research (Grants 16-01-00015, 17-51-53180) and by
the Ministry of Education and Science of Russian Federation
(Grant 1.1660.2017/4.6). The author N.P.~Bondarenko was also supported by
the Russian Federation President Grant MK-686.2017.1.\\

\begin{center}
{\bf REFERENCES}
\end{center}

\begin{enumerate}
\item[{[1]}] Hale J. Theory of functional-differential equations.
     Springer-Verlag, New York, 1977.
\item[{[2]}] Myshkis A.D. Linear differential equations with a delay
     argument. Moscow, Nauka, 1972.
\item[{[3]}] Marchenko V.A. {\it Sturm-Liouville Operators and Their
     Applications}, Naukova Dumka, Kiev, 1977; English transl.,
     Birkh\"auser, 1986.
\item[{[4]}] Levitan B.M. {\it Inverse Sturm-Liouville Problems},
     Nauka, Moscow, 1984; Engl. transl., VNU Sci.Press, Utrecht, 1987.
\item[{[5]}] Freiling G. and Yurko V.A. {\it Inverse Sturm-Liouville
     Problems and Their Applications}, NOVA Science Publishers,
     New York, 2001.
\item[{[6]}] Yurko V.A. {\it Method of Spectral Mappings in the Inverse
     Problem Theory}, Inverse and Ill-posed Problems Series. VSP,
     Utrecht, 2002.
\item[{[7]}] Freiling G. and Yurko V.A. {\it Inverse problems for
     Sturm-Liouville differential operators with a constant delay},
     Appl. Math. Lett. 25 (2012) 1999--2004.
\item[{[8]}] Yang C-Fu. Trace and inverse problem of a discontinuous
     Sturm-Liouville operator with retarded argument.
     J. Math. Anal. Appl. 395 (2012), no.1, 30-41.
\item[{[9]}] Vladi\v{c}i\'{c}, V.; Pikula, M. An inverse problem
     for Sturm-Liouville-type differential equation with a constant
     delay. Sarajevo J. Math. 12(24), no. 1 (2016), 83--88.
\item[{[10]}] Buterin, S. A.; Pikula, M.; Yurko, V. A. Sturm-Liouville
     differential operators with deviating argument, Tamkang J. Math.
     48:1 (2017), 61--71.
\item[{[11]}] Buterin, S. A.; Yurko, V. A. An inverse spectral problem
     for Sturm-Liouville operators with a large constant delay, Anal.
     Math. Phys. (2017), DOI 10.1007/s13324-017-0176-6.
\end{enumerate}

\medskip

\noindent Natalia Bondarenko \\
1. Department of Applied Mathematics, Samara National Research University, \\
Moskovskoye Shosse 34, Samara 443086, Russia, \\
2. Faculty of Mathematics and Mechanics, Saratov State University, \\
Astrakhanskaya 83, Saratov 410012, Russia, \\
e-mail: {\it bondarenkonp@info.sgu.ru}

\medskip

\noindent Vjacheslav Yurko \\
Faculty of Mathematics and Mechanics, Saratov State University, \\
Astrakhanskaya 83, Saratov 410012, Russia, \\
e-mail: {\it yurkova@info.sgu.ru}

\end{document}